\documentclass[10pt,twoside]{article}
\usepackage{graphicx}
\usepackage{amsmath}
\usepackage{Latex-document}

\def\volumeno{I}

\title{\bf  Some Recent Transcendental Techniques
\vskip -2mm in Algebraic and Complex Geometry\thanks{Research
partially supported by a grant from the National Science
Foundation.}\vskip 6mm}

\author{Yum-Tong Siu\vspace*{-0.5cm}\thanks{Department of Mathematics,
Harvard University, Cambridge, MA 02138, USA. E-mail:
siu@math.harvard.edu.}}

\date{\vspace{-8mm}}

\begin{document}

\maketitle

\thispagestyle{first} \setcounter{page}{439}

\begin{abstract}

\vskip 3mm

This article discusses the recent transcendental techniques used
in the proofs of the following three conjectures. (1)~The
plurigenera of a compact projective algebraic manifold are
invariant under holomorphic deformation. (2)~There exists no
smooth Leviflat hypersurface in the complex projective plane.
(3)~A generic hypersurface of sufficiently high degree in the
complex projective space is hyperbolic in the sense that there is
no nonconstant holomorphic map from the complex Euclidean line to
it.

\vskip 4.5mm


\end{abstract}

\vskip 12mm

\section{Introduction} \label{section 1}\setzero
\vskip-5mm \hspace{5mm }

Since the use of function theory in the study of algebraic curves as
Riemann surfaces about two hundred years ago, transcendental methods
such as harmonic forms in Hodge theory and curvature in the theory of
Chern-Weil have been very important tools in complex algebraic geometry.
Since the nineteen sixties very powerful techniques in the estimates of
$\bar\partial$, especially $L^2$ estimates and regularity techniques,
have been extensively developed by C.~B.~Morrey, J.~J.~Kohn,
L.~H\"ormander, {\it et al}. (To avoid too lengthy a bibliography here,
we refer to
\cite{C},\cite{D},\cite{Ko},\cite{Mc},\cite{SZ},\cite{S7},\cite{S6} for
references not listed here.) During the last two decades these new
transcendental techniques have been increasingly used in complex
algebraic geometry.  The most noteworthy among them is J.~J.~Kohn's
method of multiplier ideals for $\bar\partial$ estimates \cite{Ko} which
holds the promise of applicability to general partial differential
equations and global geometry. Nadel \cite{N} introduced multiplier
ideal sheaves dual to Kohn's. A number of longstanding problems in
algebraic and complex geometry hitherto beyond the reach of known
methods have been solved by the new techniques of $\bar\partial$
estimates. On the other hand, demands of geometric applications motivate
new approaches to $\bar\partial$ estimates. We will discuss here some
recent results in the following three topics in algebraic and complex
geometry obtained by the new transcendental methods. (1)~Invariance of
plurigenera. (2)~Nonexistence of smooth Levi-flat hypersurface in ${\bf
P}_2$. (3)~Hyperbolicity of generic hypersurface of high degree in ${\bf
P}_n$. Though topic~(3) is only peripherally linked to $\bar\partial$
estimates, a long outstanding problem there is solved by some recent
transcendental techniques.

\section{Invariance of plurigenera} \label{section 2}
\setzero\vskip-5mm \hspace{5mm }

Let $\Delta_r=\{t\in{\bf C}\,\big|\,|t|<r\}$ and $\Delta=\Delta_1$.  Denote by $K_Y$ the canonical line bundle of
a complex manifold $Y$. The $m$-genus of a compact complex manifold $X$ is the complex dimension of
$\Gamma\left(X,m\,K_X\right)$. By Hodge theory the $1$-genus of a compact K\"ahler manifold is a topological
invariant and therefore is invariant under holomorphic deformation.  For the general $m$-genus there is the
following conjecture on its invariance under holomorphic deformation for a compact K\"ahler manifold.

{\bf Conjecture 2.1} {\it (on Deformational Invariance of
Plurigenera for K\"ahler Manifolds). Let $\pi:X\rightarrow\Delta$
be a holomorphic family of compact K\"ahler manifolds with fiber
$X_t$.  Then for any positive integer $m$ the complex dimension of
$\Gamma\left(X_t,m\,K_{X_t}\right)$ is independent of $t$ for
$t\in\Delta$.}

Conjecture (2.1) has been verified in \cite{S3} when $X$ is a
family of compact projective algebraic manifolds.

{\bf Theorem 2.2}\ \cite{S3}. {\it Let $\pi:X\rightarrow\Delta$ be
a holomorphic family of compact complex projective algebraic
manifolds.  Then for any integer $m\geq 1$ the complex dimension
of $\Gamma\left(X_t,m\,K_{X_t}\right)$ is independent of $t$ for
$t\in\Delta$.}

The main techniques to solve the problem were first introduced in
\cite{S1} where for technical reasons the assumption of each fiber
being of general type is added. Because of the upper
semicontinuity of $\dim_{\bf C}
\Gamma\left(X_t,m\,K_{X_t}\right)$, the conjecture is equivalent
to extending every element of $\Gamma\left(X_t,m\,K_{X_t}\right)$
to $\Gamma\left(X,m\,K_X\right)$. We can assume $t=0$. The idea of
the main techniques stemmed from the following naive motivation.
If one could write an element $s^{(m)}$ of
$\Gamma\left(X_0,m\,K_{X_0}\right)$ as a sum of terms, each of
which is the product of an element $s^{(1)}$ of
$\Gamma\left(X_0,K_{X_0}\right)$ and an element $s^{(m-1)}$ of
$\Gamma\left(X_0,\left(m-1\right)K_{X_0}\right)$, then one can
extend $s^{(m)}$ to an element of $\Gamma\left(X,m\,K_X\right)$ by
induction on $m$. Of course, in general it is clearly impossible
to so express $s^{(m)}$ as a sum of such products. However, one
could successfully implement a modified form of this naive
motivation, in which $s^{(1)}$ is only a local holomorphic section
and $s^{(m-1)}$ is an element of
$\Gamma\left(X_0,\left(m-1\right)K_{X_0}+E\right)$ instead of
$\Gamma\left(X_0,\left(m-1\right)K_{X_0}\right)$, where $E$ is a
sufficiently ample line bundle on $X$ independent of $m$. The
implementation of the modified form depends on the following two
ingredients.

{\bf Proposition 2.3} {\it (Global Generation of Multiplier Ideal
Sheaves). Let $L$ be a holomorphic line bundle over an
$n$-dimensional compact complex manifold $Y$ with a Hermitian
metric which is locally of the form $e^{-\xi}$ with $\xi$
plurisubharmonic.  Let ${\cal I}_\xi$ be the multiplier ideal
sheaf of the Hermitian metric $e^{-\xi}$ ({\it i.e.,} the sheaf
consisting of all holomorphic function germs $f$ with
$\left|f\right|^2e^{-\xi}$ locally integrable). Let $E$ be an
ample holomorphic line bundle over $Y$ such that for every point
$P$ of $Y$ there are a finite number of elements of $\Gamma(Y,E)$
which all vanish to order at least $n+1$ at $P$ and which do not
simultaneously vanish outside $P$. Then $\Gamma(Y,{\cal
I}_\xi\otimes(L+E+K_Y))$ generates ${\cal I}_\xi\otimes(L+E+K_Y)$
at every point of $Y$.}

{\bf Proposition 2.4} {\it (Extension Theorem of Ohsawa-Takegoshi
Type). Let $\gamma:Y\rightarrow\Delta$ be a projective algebraic
family of compact complex manifolds. Let $Y_0=\gamma^{-1}(0)$ and
let $n$ be the complex dimension of $Y_0$. Let $L$ be a
holomorphic line bundle with a Hermitian metric $e^{-\chi}$ with
$\chi$ plurisubharmonic. Then for $0<r<1$ there exists a positive
constant $A_r$ with the following property.  For any holomorphic
$L$-valued $n$-form $f$ on $Y_0$ with $
\int_{Y_0}|f|^2e^{-\chi}<\infty$, there exists a holomorphic
$L$-valued $(n+1)$-form $\tilde f$ on $\gamma^{-1}(\Delta_r)$ such
that $\tilde f|_{Y_0}=f\wedge\gamma^*(dt)$ at points of $Y_0$ and
$ \int_Y|\tilde f|^2e^{-\chi}\leq A_r \int_{Y_0}|f|^2e^{-\chi}$.}

Locally expressing an element $s^{(m)}$ of
$\Gamma\left(X_0,m\,K_{X_0}\right)$ as a sum of terms, each of
which is the product of a local holomorphic function $s^{(1)}$ and
an element $s^{(m-1)}$ of
$\Gamma\left(X_0,\left(m-1\right)K_{X_0}+E\right)$ is precisely
Proposition 2.3, necessitating the use of $E$.

One constructs a metric of $(m-1)K_X+E$ by using the sum of
absolute-value squares of elements of
$\Gamma\left(X,\left(m-1\right)K_X+E\right)$ whose restrictions to
$X_0$ form a basis of
$\Gamma\left(X_0,\left(m-1\right)K_{X_0}+E\right)$. Proposition
2.4 is now applicable to show the surjectivity of
$\Gamma\left(X,\left(m-1\right)K_X+E\right)\to
\Gamma\left(X_0,\left(m-1\right)K_{X_0}+E\right)$ by induction on
$m$. To get rid of $E$, for a sufficiently large $\ell$ one takes
the $\ell$-th power of an element of
$\Gamma\left(X_0,m\,K_{X_0}\right)$ and multiplies it by an
element of $\Gamma\left(X,E\right)$ and then takes the $\ell$-th
root of the absolute value after its extension.  This process,
together with H\"older's inequality, is used to produce a metric
of $(m-1)K_X$ which we can use in the application of Proposition
2.4 to get the surjectivity of
$\Gamma\left(X,m\,K_X\right)\to\Gamma\left(X_0,m\,K_{X_0}\right)$.
The assumption of general type facilitates the last technical step
of getting rid of $E$ by writing $aK_X=E+D$ for some sufficiently
large integer $a$ and an effective divisor $D$.

Kawamata \cite{Ka} translated the argument of \cite{S1} to a purely
algebraic geometric setting and Nakayama \cite{Na} explored
generalizations including results concerning
$\lim_{m\to\infty}\frac{1}{m}\log \dim_{\mathbf
C}\Gamma\left(X_t,m\,K_{X_t}+E\right)$ as a function of $t$. The case of
non general type necessitates letting $\ell$, which is used in taking
the power and the root, go to infinity.  One has to control the
estimates in the limiting process.

Tsuji put on the web a preprint on the deformational invariance of
the plurigenera for manifolds not necessarily of general type
\cite{T}, in which, besides the techniques of \cite{S1}, he uses
his theory of analytic Zariski decomposition and generalized
Bergman kernels. Tsuji's approach of generalized Bergman kernels
naturally and elegantly reduces the problem of the deformational
invariance of the plurigenera to a growth estimate on the
generalized Bergman kernels.  Unfortunately this crucial estimate
is lacking and seems unlikely to be establishable, as explained in
\cite{S3}.

In \cite{S3} a metric as singular as possible is introduced for
the limiting process, which, together with an estimation technique
using the concavity of the logarithmic function, successfully
removes the technical assumption of general type in \cite{S1}.

The deformational invariance of the plurigenera for K\"ahler
manifolds is still open. Only known results on the K\"ahler case
are due to Levine's \cite{Le} with the assumption of some
pluricanonical section with nonsingular divisor (or only mild
singularities).  To generalize the methods of \cite{S1} and
\cite{S3} to the K\"ahler case, one possibility is to use the
``absolute value'' of a holomorphic line bundle constructed from
the K\"ahler metric, because in the key argument only the absolute
value of the constructed holomorphic section is used and not the
section itself.  There is still no method of implementing this
possibility.

\section{Nonexistence of smooth Levi-Flat
hypersurface in ${\bf P}_{\bf 2}$} \label{section 3}
\setzero\vskip-5mm \hspace{5mm }

The problem of the nonexistence of smooth Levi-flat hypersurface
in ${\bf P}_2$ has its origin in dynamical systems in ${\bf P}_2$
(see \cite{Li}). In terms of $\bar\partial$ estimates, its
significance is that it gives a natural geometric setting for the
understanding of $\bar\partial$ regularity for domains with
Levi-flat boundary. The $\bar\partial$ regularity problem for a
relatively compact domain $\Omega$ with smooth boundary
$\partial\Omega$ in a complex manifold is to find a solution $u$
on $\Omega$, smooth up to $\partial\Omega$, to the equation
$\bar\partial u=g$ with a given $\bar\partial$-closed $(0,1)$-form
$g$ on $\Omega$, smooth up to of $\partial\Omega$.  {\it Global
regularity} is said to hold for $\Omega$ if regularity holds for
the particular solution $u$, known as the {\it Kohn solution},
which is orthogonal to all the $L^2$ holomorphic functions on
$\Omega$. The problem of global regularity has been very
extensively studied in the past couple of decades (see
bibliographies in \cite{C},\cite{Ko}). Global regularity holds for
strictly pseudoconvex domains and, more generally, for weakly
pseudoconvex domains whose boundary points are all of finite type.
Finite type means that local complex-analytic curves touch the
boundary only to bounded finite (normalized) order. Global
regularity holds also for weakly pseudoconvex domains defined by
global smooth weakly plurisubharmonic functions. On the other
hand, worm domains are counter-examples for global regularity for
general weakly pseudoconvex domains \cite{C}. Though the
nonexistence of smooth Levi-flat hypersurface in ${\bf P}_2$ is
connected with the regularity of any one solution of the
$\bar\partial$-equation rather than the particular Kohn solution,
its proof ushers in a new approach of using vector fields to
obtain $\bar\partial$ regularity for domains with Levi-flat
boundary. The following solution of the Levi-flat hypersurface
problem was given in \cite{S4}.

{\bf Theorem 3.1} \cite{S4}. {\it Let $q\geq 8$. Then there exists
no $C^q$ Levi-flat real hypersurface $M$ in ${\bf P}_2$.}

The nonexistence of real-analytic Levi-flat hypersurface in ${\bf
P}_3$ was proved by Lins-Neto \cite{Li}. Ohsawa \cite{O} treated
the nonexistence of real-analytic Levi-flat hypersurface in ${\bf
P}_2$ (some points in the argument there not yet complete). The
nonexistence of smooth Levi-flat hypersurface in ${\bf P}_3$ was
proved in \cite{S2}. The real-analytic case is completely
different in nature from the smooth case, because the structure is
automatically extendible to a neighborhood for the real-analytic
case. Nonexistence in ${\bf P}_2$ implies nonexistence in ${\bf
P}_n$ ($n\geq 2$) by slicing with a generic linear ${\bf P}_2$.

The following argument reduces the problem to a $\bar\partial$
regularity question.  Suppose $M$ exists.  We seek a contradiction
from the positivity of the $(1,0)$-normal bundle $N_{M,{\bf
P}_2}^{(1,0)}$ of the Levi-flat hypersurface $M$.   The curvature
$\theta$ of $N_{M,{\bf P}_2}^{(1,0)}$ with the metric induced from
the Fubini-Study metric is positive, because a quotient bundle
cannot be less positive.  On the other hand, $M$ is the zero-set
of a smooth ${\bf R}$-valued function $f_M$ on ${\bf P}_2$ with
$df_M$ nowhere zero on $M$. Evaluation by $\partial f_M$ shows
that $N^{(1,0)}_{M,{\bf P}_2}$ is smoothly trivial and $\theta$
must be $d$-exact on $M$, which means that $\theta=d\alpha$ for
some smooth real $1$-form $\alpha$ on $M$. Decompose
$\alpha=\alpha^{(1,0)}+\alpha^{(0,1)}$ into its $(1,0)$ and
$(0,1)$ components. If $\alpha^{(0,1)}=\bar\partial_b\psi$ for
some smooth function $\psi$ on $M$, then
$\theta=\sqrt{-1}\partial_b\bar\partial_b\left(2\hbox{Im\,}\psi\right)$.
At a point of $M$ where $\hbox{Im\,}\psi$ assumes its maximum, the
positivity of $\theta$ along the holomorphic foliation is
contradicted.  The problem is thus reduced to solving the
$\bar\partial_b$ equation on $M$ with regularity. By applying the
Mayer-Vietoris sequence to ${\bf P}_2-M=U_1\cup U_2$ and using the
vanishing of $H^2({\bf P}_2,{\cal O}_{{\bf P}_2})$, the problem is
reduced to whether, for any $\bar\partial$-closed $(0,1)$-form $g$
on $U_j$ smooth up to $\partial U_j$, the equation $\bar\partial
u=g$ can be solved on $U_j$ with $u$ smooth up to $\partial U_j$.

The usual approach to $\bar\partial$ regularity is to use the
Bochner-Kodaira formula with boundary $\|\bar\partial
g\|_\Omega^2+\|\bar\partial^*g\|_\Omega^2
=\int_{\partial\Omega}\left<{\cal L},\bar g\wedge
g\right>+\|\bar\nabla g\|_\Omega^2 +\left(\Theta_E,\bar g\wedge
g\right)_\Omega$ (with $g$ being a smooth $E$-valued $(n,1)$-form
in the domain of $\bar\partial^*$), to solve the equation with
$L^2$ estimates and then apply differential operators, integration
by parts, and commutation relations to prove regularity. Here
$n=\dim_{\bf C}\Omega$, $\|\cdot\|_\Omega$ is the $L^2$ norm over
$\Omega$, $\bar\partial^*$ is the adjoint of $\bar\partial$,
$\bar\nabla$ means covariant differentiation in the
$(0,1)$-direction, ${\cal L}$ is the Levi form of
$\partial\Omega$, and $\Theta_E$ is the curvature form of the
Hermitian holomorphic line bundle $E$ (which is usually chosen to
be trivial).

In our new approach to get $\bar\partial$ regularity for the
Levi-flat domain $\Omega$ in ${\bf P}_2$ the use of holomorphic
vector fields compensates for the complete lack of strict
positivity for the Levi form of the boundary.  We use a new norm
to derive the Bochner-Kodaira formula with boundary.  We choose a
vector field $\xi$ on ${\bf P}_2$ which generates biholomorphisms
preserving the Fubini-Study metric. The new norm is the $L^2$ norm
$L_m^2(\Omega,\xi)$ for Lie derivatives $\left({{\cal
L}ie}_\xi\right)^j g$ along $\xi$ for order $j\leq m$ on $\Omega$
for $(0,1)$-form $g$.  Since $\xi$ generates metric-preserving
biholomorphisms, the formal adjoint of $\bar\partial$ with respect
to $L_m^2(\Omega,\xi)$ agrees with the one with respect to usual
$L^2$.  One usual difficulty with regularity is the error terms
from the commutation of differential operators with $\bar\partial$
and $\bar\partial^*$.  One advantage of using ${{\cal L}ie}_\xi$
is that there are no error terms from its commutation with
$\bar\partial$ and $\bar\partial^*$.

There are two technical problems.  One is how to establish, for
such a norm, the Bochner-Kodaira formula with boundary. The other
is that appropriate regularity for a solution of the
$\bar\partial$ equation with finite $L_m^2(\Omega,\xi)$ norm can
be obtained only at points where the real and imaginary parts of
$\xi$ are not both tangential to $\partial\Omega$.  We handle the
first problem as follows. We prove that, if $g$ belongs to the
domain of the adjoint of $\bar\partial$ with respect to
$L_m^2(\Omega,\xi)$, then $\left({{\cal L}ie}_\xi\right)^j g$
belongs to the domain of the adjoint of $\bar\partial$ with
respect to the usual $L^2$ norm on $\Omega$ for $j\leq m$. The
formula for the new norm is simply the sum, over $0\leq j\leq m$,
of such a formula for the usual $L^2$ norm on $\Omega$ for
$\left({{\cal L}ie}_\xi\right)^j g$. The proof for $\left({{\cal
L}ie}_\xi\right)^j g$ to belong to the domain of the adjoint of
$\bar\partial$ with respect to the usual $L^2$ norm on $\Omega$
consists of two steps. One shows that this is locally true at
points where $\xi$ is not tangential to $\partial\Omega$. Then one
uses a removable singularity argument to handle the other points
when $\xi$ has been chosen generic enough. For the second problem,
to handle the other points for a generic $\xi$, we use the
foliation of $\partial\Omega$ by local complex-analytic curves and
the generalized Cauchy integral formula along the local
complex-analytic curves.

\section{Hyperbolicity of generic hypersurface
of high degree in {\boldmath ${\rm P}_n$}} \label{section 4}
\setzero\vskip-5mm \hspace{5mm }

A complex manifold $X$ is {\it hyperbolic} if there
exists no nonconstant holomorphic map ${\bf C}\rightarrow X$. For
the last few decades the study of hyperbolicity has been focussed
on hypersurfaces and their complements in two important settings:
(1)~inside an abelian variety and (2)~inside ${\bf P}_n$. In the
general setting hyperbolicity is conjectured to be linked to the
positivity of canonical line bundle in the following formulation.

{\bf Conjecture 4.1}{\it (Conjecture of Green-Griffiths). In a
compact algebraic manifold $X$ of general type (or with positive
canonical line bundle) there exists a proper subvariety $Y$
containing the images of all nonconstant holomorphic maps ${\bf
C}\rightarrow X$.}

The theory for the setting inside an abelian variety is very well
developed (see \cite{S7},\cite{S6} for references). The Zariski closure
of any holomorphic map from ${\bf C}$ to an abelian variety $A$ is the
translate of an abelian subvariety of $A$. In particular, a subvariety
of an abelian variety $A$ which does not contain any translate of an
abelian subvariety of $A$ is hyperbolic. The defect of an ample divisor
in an abelian variety is zero. In particular, the complement of an ample
divisor in an abelian variety is hyperbolic.

Except those motivated by methods of number theory due to McQuillan,
practically all the major techniques for problems related to
hyperbolicity in the setting of abelian varieties are due to Bloch
\cite{B} who introduced the use of holomorphic jet differentials and
differential equations in conjunction with the jet differentials.
Investigations on problems related to hyperbolicity in the setting of
abelian varieties have essentially been completed. Only technical
details such as getting an optimal lower bound for $k_n$ in Theorem 4.2
below remain open. Theorem 4.2 (proved in Addendum of \cite{SY}) was
added to \cite{SY} in response to a difficulty in the proof of Lemma 2
of the original paper \cite{SY} pointed out in \cite{No}. The difficulty
resulted from an attempt to use semi-continuity of cohomology groups in
deformations to avoid employing Bloch's technique from \cite{B} which
involves the uniqueness part of the fundamental theorem of ordinary
differential equations.  Putting back Bloch's technique removes the
difficulty and at the same time improves the zero defect statement in
\cite{SY} to Theorem 4.2 on the second main theorem with truncated
multiplicity.

{\bf Theorem 4.2} (Addendum, \cite{SY}). {\it Let $D$ be an ample
divisor of an abelian variety $A$ of complex dimension $n$ and let
$k_0=0$, $k_1=1$, and
$k_{\ell+1}=k_\ell+3^{n-\ell-1}\left(4(k+1)\right)^\ell D^n$\ (
$1\leq\ell<n$). Then for any holomorphic map $f:{\bf C}\rightarrow
A$ whose image is not contained in any translate of $D$, the
following second main theorem with truncated multiplicity holds:
$m(r,f,D)+\left(N(r,f,D)-N_{k_n}(r,f,D)\right) =O(\log
T(r,f,D)+\log r)$ for $r$ outside some set whose measure with
respect to $\frac{dr}{r}$ is finite.  Here $T(r,f,D)$, $m(r,f,D)$,
$N(r,f,D)$, $N_{k_n}(r,f,D)$ are respectively the characteristic,
proximity, counting functions, and truncated counting functions.}

For the setting inside ${\bf P}_n$ there is the following outstanding
conjecture.

{\bf Conjecture 4.3} (Kobayashi's Conjecture). {\it {\rm (a)} The
complement in ${\bf P}_n$ of a generic hypersurface of degree at
least $2n+1$ is hyperbolic. {\rm (b)} A generic hypersurface of
degree at least $2n-1$ in ${\bf P}_n$ is hyperbolic for $n\geq
3$.}

For Conjecture 4.3(a) the complement in ${\bf P}_2$ of a generic curve
of sufficiently high degree is known to be hyperbolic \cite{SY1}. For
Conjecture 4.3(b) a generic surface of degree $\geq 36$ in ${\bf P}_3$
is known to be hyperbolic \cite{Mc}.  The degree bound is lowered to
$21$ in \cite{D}. There are some constructions of examples of smooth
hyperbolic hypersurfaces in ${\bf P}_n$ (see \cite{SZ}).  The
hyperbolicity result we want to discuss here is the following.

{\bf Theorem 4.4} \cite{S5}. {\it There exists a positive integer
$\delta_n$ such that a generic hypersurface in ${\bf P}_n$ of
degree $\geq\delta_n$ is hyperbolic.}

We sketch its proof here. A central role will be played by jet
differentials which we now define. A $k$-jet differential on a
complex manifold $X$ with local coordinates $x_1,\cdots,x_n$ is
locally a polynomial in $d^\ell x_j\ (1\leq\ell\leq k,1\leq j\leq
n)$.

{\bf Lemma 4.5} (Lemma of Jet Differentials). {\it If a
holomorphic jet differential $\omega$ on a compact complex
manifold $X$ vanishes on an ample divisor of $X$ and $\varphi:{\bf
C}\rightarrow X$ is a holomorphic map, then $\varphi^*\omega$ is
identically zero on ${\bf C}$.}

The intuitive reason for Lemma 4.5 is that ${\bf C}$ does not admit a
metric (or not even a {\it $k$-jet metric}) with curvature bounded above
by negative number. While a usual metric assigns a value to a tangent
vector (which is a $1$-jet), a $k$-jet metric assigns a value to a
$k$-jet.  A non identically zero $\varphi^*\omega$ defines a $k$-jet
metric $\left|\varphi^*\omega\right|^2$ on ${\bf C}$ which, even with
some degeneracy, still gives a contradiction by its negative curvature.
A rigorous proof of Lemma 4.5 depends on the logarithmic derivative
lemma of Nevanlinna theory.  A consequence of Lemma 4.5 is that the
image of the $k$-jet $d^k\varphi$ of any holomorphic map $\varphi:{\bf
C}\rightarrow X$ satisfies the differential equation $\omega=0$ on $X$.
If there exist enough independent such $\omega$ on $X$, then the system
of all equations $\omega=0$ does not admit any local solution curve and
$X$ is hyperbolic.

In the setting of abelian varieties Bloch constructed jet differentials
by comparing meromorphic functions on the image and the target of a map
with finite fibers. For a holomorphic map $\varphi$ from ${\bf C}$ to an
abelian variety $A$, let $X$ be the Zariski closure of the image of
$\varphi$ in $A$ and $\cal X$ be the Zariski closure of
$(d^k\varphi)({\bf C})$ in $J_k(A)=A\times{\bf C}^{nk}$.  Here
$J_k(\cdot)$ means the space of $k$-jets.  Let $\sigma_k:{\cal X}\to
{\bf C}^{nk}$ be induced by the natural projection $J_k(A)=A\times{\bf
C}^{nk}\rightarrow{\bf C}^{nk}$ which forgets the position and keeps the
differentials. Let $\tau: J_k(X)\rightarrow X$ be the natural
projection. Let $F$ be a meromorphic function on $X$ whose pole-set is
some ample divisor $D$. Suppose $\sigma_k:{\cal X}\rightarrow{\bf
C}^{nk}$ is generically finite. Let $x_1,\cdots,x_n$ be the coordinates
of ${\bf C}^n$. Then $F\circ\tau$ belongs to a finite extension of the
rational function field of ${\bf C}^{nk}$ and there exist polynomials
$P_j$ ($0\leq j\leq p$) with constant coefficients in the variables
$d^\ell x_\nu$ ($1\leq\ell\leq k,1\leq\nu\leq n$) such that $
\sum_{j=0}^p (\sigma_k^*P_j)(\tau^*F)^j=0$ on ${\cal X}$ and
$\sigma_k^*P_p$ is not identically zero on ${\cal X}$. The equation
forces the holomorphic jet differential $P_p$ on ${\cal X}$ to vanish on
the ample divisor $\tau^{-1}(D)$. The assumption of generical finiteness
of $\sigma_k$ is tied to the translational invariance of $X$.

The idea of our method of construction of holomorphic jet differentials
on a generic hypersurface $X$ in ${\bf P}_n$ defined of by a polynomial
$f$ of degree $\delta$ is to use the theorem of Riemann-Roch and the
lower bound of negativity of jet differential bundles of $X$.  The
theorem of Riemann-Roch was first used by Green-Griffiths to obtain
holomorphic jet differentials and is applicable only for surfaces where
the higher cohomology groups could be easily handled.  We can handle the
higher cohomology groups in our higher dimensional case because of the
lower bound of negativity of jet differential bundles of $X$.  Since the
twisted cohomology groups of ${\bf P}_n$ comes from counting the number
of monomials, in the actual proof direct counting of monomials is used.
Let $x_1,\cdots,x_n$ (respectively $z_0,\cdots,z_n$) be the
inhomogeneous (respectively homogeneous) coordinates of ${\bf P}_n$.
Let $Q$ be a non identically zero polynomial of degree $m_0$ in
$x_1,\cdots,x_n$ and of homogeneous weight $m$ in $d^jx_\ell$ ($1\leq
j\leq n-1,\,1\leq\ell\leq n$) with the weight of $d^jx_\ell$ equal to
$j$. If $m_0+2m<\delta$, then $Q$ is not identically zero on
$J_{n-1}(X)$. By counting the number of coefficients of $Q$ and the
number of equations needed for the jet differential on $X$ defined by
$Q$ to vanish on an ample divisor in $X$ of high degree defined by a
polynomial $g=0$ in ${\bf P}_n$ and using $f=df=\cdots=d^{n-1}f=0$ to
eliminate one coordinate and its differentials, we obtain a jet
differential $\frac{Q}{g}$ on $X$ which is holomorphic and vanishes on
an ample divisor of high degree.

{\bf Proposition 4.6} (Existence of Holomorphic Jet
Differentials). {\it If
$0<\theta_0,\,\theta,\,\theta^\prime<1-\epsilon$ with
$n\theta_0+\theta\geq n+\epsilon$, then there exists an explicit
$A=A(n,\epsilon)>0$ such that for $\delta\geq A$ there exists a
non identically zero ${\cal O}_{{\bf P}_n}(-q)$-valued holomorphic
$(n-1)$-jet differential $\omega$ on $X$ of total weight $m$ with
$q\geq\delta^{\theta^\prime}$ and $m\leq\delta^\theta$.}

To construct enough independent jet differentials, we use meromorphic
vector fields of low pole order on the total space ${\cal X}$ of all
hypersurfaces in ${\bf P}_n$ of degree $\delta$. The total space ${\cal
X}$ is defined by $f=\sum_{\nu\in{\bf N}^{n+1},|\nu|=\delta} \alpha_\nu
z^\nu$ of bidegree $(\delta,1)$ in ${\bf P}_n\times{\bf P}_N$, where
$N={\delta+n\choose n}-1$, $z^\nu=z_0^{\nu_0}\cdots z_n^{\nu_n}$,
$|\nu|=\nu_0+\cdots+\nu_n$, and ${\bf N}$ is the set of all nonnegative
integers. Let $e_\ell=(0,\cdots,0,1,\cdots,0)\in{\bf N}^{n+1}$ with $1$
in the $\ell$-th place.  The $(1,0)$-twisted tangent bundle of ${\cal
X}$ is globally generated by holomorphic sections of the forms
$L\left(z_q\left(\frac{\partial}{\partial\alpha_{\lambda+e_p}}\right)
-z_p\left(\frac{\partial}{\partial\alpha_{\lambda+e_q}}\right)\right) $
and $\sum_j B_j\frac{\partial}{\partial z_j}+\sum_\mu
L_\mu\frac{\partial}{\partial\alpha_\mu}$, where $\lambda\in{\bf
N}^{n+1}$ with $|\lambda|=\delta-1$ and $L,\,L_\mu$ (respectively $B_j$)
are homogeneous linear functions of $\{\alpha_\nu\}$ (respectively
$z_0,\cdots,z_n$) with $L_\mu$ and $B_j$ suitably chosen.

We introduce the space $J^{\rm{\scriptstyle vert}}_{n-1}\left({\cal
X}\right)$ of vertical $(n-1)$-jets of ${\cal X}$ which is defined by $
f=df=\cdots=d^{n-1}f=0 $ in $\left(J_{n-1}\left({\bf
P}_n\right)\right)\times{\bf P}_N$ with the coefficients $\alpha_\nu$ of
$f$ regarded as constants when forming $d^jf$. By generalizing the above
construction of vector fields on ${\cal X}$ to vector fields on
$J^{\rm{\scriptstyle vert}}_{n-1}\left({\cal X}\right)$, one obtains the
following.

{\bf Proposition 4.7} (Existence of Low Pole-Order Vector Fields).
{\it There exist $c_n,\, c^\prime_n\in{\bf N}$ such that the
$(c_n,c_n^\prime)$-twisted tangent bundle of the projectivization
of $J^{\rm{\scriptstyle vert}}_{n-1}\left({\cal X}\right)$ is
globally generated. (To avoid considering the singularities of
weighted projective spaces, one can interpret the statement by
using functions which are polynomials of homogeneous weight along
the fibers of $J^{\rm{\scriptstyle vert}}_{n-1}\left({\cal
X}\right)$.)}

For a generic fiber $X$ of ${\cal X}$ the constructed holomorphic
$(n-1)$-jet differential $\omega$ on $X$ with vanishing order at least
$q$ on the infinity divisor can be extended holomorphically to
$\tilde\omega$ on all neighboring fibers with vanishing order at least
$q$ on the infinity divisor. We use vector fields $v_1,\cdots,v_p$ on
$J^{\scriptstyle{\rm vert}}_{n-1}\left({\cal X}\right)$ with fiber pole
order low relative to $q$ and take successive Lie derivatives ${\cal
L}ie_{v_1}\cdots{\cal L}ie_{v_p}\tilde\omega$ whose restrictions to $X$
give holomorphic jet differentials on $X$ vanishing on an ample divisor
of $X$. Because of the bound on the weight $m$ in the construction of
$\omega$, for $\delta$ sufficiently large the jet differentials from the
Lie derivatives are independent enough to eliminate the derivatives from
the differential equations they define.  As a result, one concludes that
for some proper subvariety $Y$ in $X$ the image of any nonconstant
holomorphic map from ${\bf C}$ to $X$ is contained in $Y$.

To get the full conclusion of hyperbolicity, for the constructed
$\omega$ one has to control the vanishing order of the coefficients of
the monomials of the differentials.  For a generic $X$ the construction
process enables one to bound the vanishing order by $\delta^{2-\eta}$
for some $\eta>0$.  For hyperbolicity one needs the better bound of
$\delta^{1-\eta}$ for some $\eta>0$. To achieve it, one uses an
appropriate embedding $\Phi:{\bf P}_n\to{\bf P}_{\hat n}$ of degree
$\delta_1$ so that a generic hypersurface $X$ of degree
$\delta:=\delta_1\delta_2$ in ${\bf P}_n$ can be extended to a
hypersurface $\hat X$ of degree $\delta_2$ in ${\bf P}_{\hat n}$. For
this step the method of multiplier ideal sheaves from $\bar\partial$
estimates is used. We deform $\Phi$ slightly and pull back the jet
differential $\hat\omega$ constructed on $\hat X$ to get a differential
$\omega$ on a slight deformation of $X$. When the image of the deformed
$\Phi$ has appropriate transversality to the zero set of the
coefficients of $\hat\omega$, an appropriate choice of $\delta_1$ and
$\delta_2$ gives the required bound on the vanishing order of the
coefficients of $\omega$.  This is at the expense of increasing the
order of the jet differential from $n-1$ to $\hat n-1$, which does not
affect the argument. For this additional step in the argument the degree
$\delta$ must be a product. To remove this condition, one uses an
embedding ${\bf P}_n\to{\bf P}_{\hat n_1}\times{\bf P}_{\hat n_2}$
instead of $\Phi$.

The use and the construction of meromorphic vector fields on
$J^{\scriptstyle{\rm vert}}_{n-1}\left({\cal X}\right)$ of low pole
order along the fibers are motivated by Clemens's work \cite{Cl} (with
later generalizations and improvements by Ein \cite{E} and Voisin
\cite{V}) on the nonexistence of regular rational and elliptic curves on
generic hypersufaces of sufficiently high degree.

There is no way yet to handle Conjecture 4.1.  Additional assumptions
such as $K_X-mL>0$ or $\left(K_X-mL\right)L^{n-1}>0$ for some large $m$
and $L$ ample or very ample could facilitate the construction of
holomorphic jet differentials vanishing on an ample divisor. One
possibility to handle the question of enough independent jet
differentials is to deform $d^k\varphi$ of $\varphi$ for each $k\geq 1$
separately and use techniques analogous to the twisted difference maps
in the Vojta-Faltings proof of the Mordell conjecture and to McQuillan's
separate rescaling of an entire holomorphic curve in each factor of a
product of several copies of an abelian variety (see
pp.504-505,\cite{S7}).

\label{lastpage}

\end{document}